\documentclass[12pt]{article}
\usepackage{color} %
\usepackage{amssymb,amsfonts} %

\def\proofr{\begin{pf}}
\def\proofend{$\Box$\end{pf}}

\usepackage{graphicx} %

\renewenvironment{abstract}{\par\indent\small \textbf{Abstract.}}{\par\normalsize}
\newenvironment{frontmatter}{}{}
\newenvironment{keyword}{\emph{Keywords}:}{}
\providecommand{\sep}{, }
\providecommand{\MSC}{\par\emph{MSC}:\ }
\renewcommand{\title}[1]{\begin{center}{\Large\textbf{#1}}\\[2.0ex]}
\renewcommand{\author}[1]{{\Large\rm #1}\end{center}}
\providecommand{\ead}[1]{}

\providecommand{\address}[1]{\footnotetext{Denis Krotov, #1. \texttt{krotov@math.nsc.ru}}}
\providecommand{\journal}[1]{}
\newenvironment{Proof}[1][\hspace{-1ex}]%
 {\par\addvspace{0.3em}{\bf Proof\hspace{1ex}#1.~~}}%
 {~~$\Box$\par\addvspace{0.3em}}

\newcounter{thrm}
\newcounter{lmm}

\newcounter{crll}
\newcounter{rmrk}
\newenvironment{Theorem}[1][\hspace{-1ex}]%
 {\par\addvspace{0.3em}\refstepcounter{thrm}\bf Theorem~\thethrm\hspace{1ex}{\rm #1}.~~\it}%
 {\rm\par\addvspace{0.3em}}

\newenvironment{Lemma}[1][\hspace{-1ex}]%
 {\par\addvspace{0.3em}\refstepcounter{lmm}\bf Lemma~\thelmm\hspace{1ex}{\rm #1}.~~\it}%
 {\rm\par\addvspace{0.3em}}

\newenvironment{Corollary}[1][\hspace{-1ex}]%
 {\par\addvspace{0.3em}\refstepcounter{crll}\bf Corollary~\thecrll\hspace{1ex}{\rm #1}.~~\it}%
 {\rm\par\addvspace{0.3em}}

\newenvironment{Note}[1][\hspace{-1ex}]%
 {\par\addvspace{0.3em}\refstepcounter{rmrk}\bf Remark~\thermrk\hspace{1ex}{\rm #1}.~~\rm}%
 {\rm\par\addvspace{0.3em}}

\journal{European Journal of Combinatorics}

\def\eqdf{\stackrel{\rm\scriptscriptstyle def}=}
\def\iffdf{\,\stackrel{\rm\scriptscriptstyle def}\Longleftrightarrow\,}
\def\ifff{\,\iff\,}
\begin{document}
\begin{frontmatter}
\title{On irreducible $n$-ary quasigroups with reducible 
retracts}
\author{Denis Krotov}
\ead{krotov@math.nsc.ru}
\address{Sobolev Institute of Mathematics, pr-t Ak. Koptyuga, 4,
Novosibirsk, 630090, Russia}
\begin{abstract}
An $n$-ary operation $Q:\Sigma^n\to \Sigma$ is called
an {$n$-ary quasigroup} {of order $|\Sigma|$} if
in $x_{0}=Q(x_1, \ldots , x_n)$ knowledge of any $n$ elements
of $x_0$, \ldots , $x_n$ uniquely specifies the remaining one.
An $n$-ary quasigroup $Q$ is permutably reducible if
$Q(x_1,\ldots,x_n)=P\left(R(x_{\sigma (1)},\ldots,x_{\sigma (k)}),
\linebreak[1]
x_{\sigma (k+1)},\ldots,x_{\sigma (n)}\right)$
where $P$ and $R$ are $(n-k+1)$-ary and $k$-ary quasigroups, $\sigma$ is a permutation,
and $1<k<n$.
For even $n$ we construct a permutably irreducible $n$-ary quasigroup of order $4r$
such that all its retracts obtained by fixing one variable are permutably reducible.
We use a partial Boolean function that satisfies similar
properties. For odd $n$ the existence of a permutably irreducible $n$-ary quasigroup
such that all its $(n-1)$-ary retracts are permutably reducible is an open question;
however, there are nonexistence results
for $5$-ary and $7$-ary quasigroups of order~$4$.

\end{abstract}
\end{frontmatter}


\section{Introduction}\label{s:intro}

An $n$-ary operation $Q:\Sigma^n\to \Sigma$,
where $\Sigma$ is a nonempty set,
is called
an \emph{$n$-ary qua\-si\-group} or \emph{$n$-qua\-si\-group}
(\emph{of order $|\Sigma|$}) if 
in the equality $z_{0}=Q(z_1, \ldots , z_n)$ knowledge of any $n$ elements
of $z_0$, $z_1$, \ldots , $z_n$ uniquely specifies the remaining one \cite{Belousov}.
The definition is symmetric with respect to the variables
$z_0$, $z_1$, \ldots , $z_n$, and sometimes it is comfortable to use
a symmetric form for the equation $z_{0}=Q(z_1, \ldots , z_n)$.
For this reason, we will write
\begin{equation}\label{e:[]}
Q\langle z_0, z_1, \ldots , z_n\rangle \iffdf z_0 = Q(z_1, \ldots , z_n).
\end{equation}
If we assign some fixed values to $l\leq n$ variables
in the predicate $Q\langle z_0, \ldots , z_n\rangle$
then the $(n-l+1)$-ary predicate obtained corresponds to an $(n-l)$-qua\-si\-group.
Such a qua\-si\-group is called a \emph{retract} of $Q$.
We say that an $n$-qua\-si\-group $Q$ is \emph{$A$-re\-du\-ci\-ble} if 
\begin{equation}\label{e:red}
Q \langle z_0,\ldots, z_{n} \rangle
\ifff
Q'(z_{a_1},\ldots,z_{a_k})=Q''(z_{b_1},\ldots,z_{b_{n-k+1}})
\end{equation}
where $A=\{a_1, \ldots ,a_k\}=\{0, \ldots ,n\} \backslash \{b_1,  \ldots ,b_{n-k+1}\}$
and $Q'$ and $Q''$ are $k$- and $(n-k+1)$-qua\-si\-groups.
An $n$-qua\-si\-group is \emph{permutably reducible} if 
it is $A$-re\-du\-ci\-ble for some
$A\subset \{0, \ldots ,n\}$, $1<|A|<n$.
In what follows we omit the word ``permutably'' because we consider only that type of reducibility
(often, ``reducibility'' of $n$-qua\-si\-groups denotes
the so-called $(i,j)$-re\-du\-ci\-bi\-li\-ty,
see Remark~\ref{r:red}).
In other words, an $n$-qua\-si\-group is re\-du\-ci\-ble if 
it can be represented as a repetition-free superposition of qua\-si\-groups with smaller arities.
An $n$-qua\-si\-group is \emph{irreducible} if 
it is not re\-du\-ci\-ble.

In \cite{Kro:2004ACCT:decomp,Kro:n-3}, it was shown that if the maximum
arity $m$ of an irreducible retract of an $n$-qua\-si\-group $Q$ belongs to $\{3, \ldots , n-3\}$
then $Q$ is reducible. Nevertheless, this interval does not contain $2$ and $n-2$
and thus can not guarantee the nonexistence of
an irreducible $n$-qua\-si\-group all of whose $(n-1)$-ary retracts are reducible.
In this paper we show that, in the case of order $4r$, such an $n$-qua\-si\-group
exists for 
even $n\geq 4$.
In the case of odd $n$,
as well as in the case of orders that are not divisible by $4$,
the question remains open;
however, as the result of an exhaustive computer search, we can state the following:
\begin{itemize}
    \item
There is no irreducible $5$- or $7$-qua\-si\-group of order $4$ such that
all its $(n-1)$-ary retracts are reducible.
\end{itemize}
For given order, constructing irreducible $n$-qua\-si\-groups
with reducible $(n-1)$-ary retracts is a more difficult task than simply
constructing irreducible
$n$-qua\-si\-groups.
In the last case we can break the reducibility of an $n$-qua\-si\-group
by changing it locally \cite{KroPot:ir}.
For our aims local modifications do not work properly because they
also break the reducibility of retracts.

In Section~\ref{s:q} we use a variant of the product of $n$-qua\-si\-groups of order $2$
to construct $n$-quasigroups of order $4$ from partial Boolean functions defined on the
even (or odd) vertices of the Boolean $(n+1)$-cu\-be.
The class constructed plays an important role for the $n$-qua\-si\-groups of order $4$;
up to equivalence, it gives almost all $n$-qua\-si\-groups of order $4$, see \cite{PotKro:asymp}.
It turns out that the reducibility of such an $n$-qua\-si\-group is equivalent
to a similar property, separability, of the corresponding partial Boolean function.
So, for this class the main question is reduced
to the same question for partial Boolean functions.
In Section~\ref{s:b} we construct a partial Boolean function with the required properties.
In Section~\ref{s:g} we consider the graph interpretation of the result.
\section{{$\lowercase{n}$}-Quasigroups of order $4$ and partial Boolean functions} \label{s:q}

In this section we consider $n$-quasigroups over the set
$\Sigma = Z_2^2 = \{[0,0],\linebreak[1] [0,1],\linebreak[1] [1,0],\linebreak[1] [1,1] \}$
and partial Boolean functions defined on the following subsets
of the Boolean hypercube $E^{n+1}\eqdf\{0,1\}^{n+1}$:
$$E^{n+1}_\alpha \eqdf \{(x_0,\ldots,x_n)\in E^{n+1} \,|\, x_0+\ldots+x_n=\alpha\},\quad \alpha\in\{0,1\}.$$
All calculations with elements of $\{0,1\}$ are made modulo~$2$,
while all calculations with indices are modulo~$n+1$, for example,
$x_{-1}$ means the same as $x_{n}$.
Note that, since any coordinate (say, the $0$th) in $E^{n+1}_0$ is the sum of the others,
partial Boolean functions defined on $E^{n+1}_0$ (as well as on $E^{n+1}_1$)
can be considered as Boolean functions on $E^{n}$; however, the form
that is symmetrical with respect to all $n+1$ coordinates helps to improve the presentation,
as in the case of $n$-qua\-si\-groups.

We will use the following notation:
if $j\geq i$ then
\begin{itemize}
    \item
$\overline {i,j}$ means $i, {i+1}, \ldots ,j$;
    \item
$x_{i}^{j}$ means $x_i, x_{i+1}, \ldots ,x_j$;
    \item
$|x_{i}^{j}|$ means the sum $x_i+x_{i+1}+ \ldots +x_j$;
    \item
$[x,y]_{i}^{j}$ means $[x_i,y_i], [x_{i+1},y_{i+1}], \ldots , [x_j,y_j]$;
    \item
$0^k$ means $k$ zeroes.
\end{itemize}
Given $\alpha\in\{0,1\}$ and $\lambda:E^{n+1}_{\alpha}\to \{0,1\}$,
define the $n$-qua\-si\-group
$Q_{\alpha,\lambda}$ as
\begin{equation}\label{e:constr}
Q_{\alpha,\lambda}\langle[x,y]_{0}^{n}\rangle
\iffdf
\cases{
|x_{0}^{n}| =\alpha,
\cr
|y_{0}^{n}| = \lambda(x_{0}^{n})
}
\end{equation}
or, equivalently,
\begin{equation}\label{e:constr2}
Q_{\alpha,\lambda}([x,y]_{1}^{n})
\eqdf
\left[\,|x_{1}^{n}| + \alpha,\ |y_{1}^{n}| + \dot\lambda(x_{1}^{n})\,\right]
\end{equation}
where $\dot\lambda(x_{1}^{n})\eqdf\lambda(|x_{1}^{n}|+\alpha,x_{1}^{n})$
is a representation of $\lambda$ as a Boolean function $E^n\to\{0,1\}$.
Note that we will use $\alpha$ only in the proof of Theorem~\ref{t:q}(b,c),
and it is not needed for formulating the main result.
In Lemma~\ref{l:0} below, we will see that the reducibility property of $Q_{\alpha,\lambda}$
corresponds to a similar property of the function $\lambda$.

We say that a partial Boolean function $\lambda:E^{n+1}_{\alpha}\to \{0,1\}$
is \emph{$A$-se\-pa\-rable} if 
\begin{equation}\label{e:sep}
\lambda(x_{0}^{n})\equiv \lambda'(x_{a_1},\ldots,x_{a_k})+\lambda''(x_{b_1},\ldots,x_{b_{m}})
\end{equation}
where $A=\{a_{1}^{k}\}=\{\overline{0,n}\} \backslash \{b_{1}^{m}\}$
and $\lambda':E^k\to \{0,1\}$, $\lambda'':E^m\to \{0,1\}$ are Boolean functions.
(Here and elsewhere $\equiv $ means that the two expressions are identical
on the region of the left one.)
$\lambda$ is \emph{separable} if 
it is $A$-se\-pa\-rable for some
$A\subset \{\overline{0,n}\}$, $2\leq|A|\leq n-1$.

\begin{Lemma}\label{l:0}
Let $A\subset \{\overline{0,n}\}$.
The $n$-qua\-si\-group $Q_{\alpha,\lambda}$ is $A$-re\-du\-ci\-ble if and only if
the partial Boolean function $\lambda:E^{n+1}_{\alpha}\to \{0,1\}$ is $A$-se\-pa\-rable.
\end{Lemma}
In the proof, we will use the following simple fact \cite{Kro:2004ACCT:decomp,Kro:n-3}:
\begin{Lemma} \label{l:abg}
Assume two $n$-qua\-si\-groups $Q_1$ and $Q_2$ are $\{\overline{0,k-1}\}$-re\-du\-ci\-ble.
If $Q_1\langle z_{0}^{k-1},~z_k,0^{n-k} \rangle \iff Q_2\langle z_{0}^{k-1},~z_k,0^{n-k}\rangle$
and $Q_1\langle z_0,0^{k-1},~z_{k}^{n} \rangle \iff Q_2\langle z_0,0^{k-1},~z_{k}^{n}\rangle$
then $Q_1$ and $Q_2$ are identical.
\end{Lemma}

\begin{Proof}[of Lemma~\ref{l:0}]
Clearly, (\ref{e:sep}) implies (\ref{e:red}) with $Q=Q_{\alpha,\lambda}$ (see (\ref{e:constr})),
and $Q'=Q_{\alpha,\mu}$, $Q''=Q_{0,\nu}$ where $\dot\mu=\lambda'$, $\dot\nu=\lambda''$ (see (\ref{e:constr2})).

Let us prove the converse.
Suppose $Q_{\alpha,\lambda}$ is $A$-reducible.
Without loss of generality assume $\alpha=0$ and $A=\{\overline{0,k-1}\}$.
Using Lemma~\ref{l:abg}, we can verify that $Q_{0,\lambda}\langle[x,y]_{0}^{n}\rangle$
defined by (\ref{e:constr}) is equivalent to
$$
\cases{
|x_{0}^{n}|=0,\cr
|y_{0}^{n}|=\lambda (x_{0}^{k-1},|x_{0}^{k-1}|,0^{n-k})
+\lambda(|x_{0}^{k-1}|,0^{k-1},|x_{0}^{k-1}|,0^{n-k})+\lambda (|x_{k}^{n}|,0^{k-1},x_{k}^{n}).}
$$
Comparing with (\ref{e:constr}), we find that
$ \lambda(x_{0}^{n})\equiv \lambda'(x_{0}^{k-1})+\lambda''(x_{k}^{n})$
where
\begin{eqnarray*}
\lambda'(x_{0}^{k-1})&\eqdf&\lambda (x_{0}^{k-1},|x_{0}^{k-1}|,0^{n-k})
+\lambda(|x_{0}^{k-1}|,0^{k-1},|x_{0}^{k-1}|,0^{n-k}),\\
\lambda''(x_{k}^{n}) &\eqdf& \lambda (|x_{k}^{n}|,0^{k-1},x_{k}^{n}).
\end{eqnarray*}
Therefore $\lambda$ is $\{\overline{0,k-1}\}$-se\-pa\-rable.
\end{Proof}
The following main theorem results from Lemma~\ref{l:0} and Theorem~\ref{t:b}
from the next section.
Although the proof 
depends on Theorem~\ref{t:b},
it is straightforward, and placing it first hardly leads to mishmash.

\begin{Theorem}\label{t:q}
Let $n\geq 4$ be even and
$f(x_{0}^{n}) \eqdf \sum_{i=0}^{n}\sum_{i=1}^{\lfloor n/4 \rfloor}x_i x_{i+j}$
for all $x_{0}^{n}\in E^{n+1}_0$.
Then\\
{\rm (a)} The $n$-qua\-si\-group $Q_{0,f}$ is irreducible.\\
{\rm (b)} Every $(n-1)$-ary retract $Q^i_{[\alpha,\gamma]}$
obtained from $Q_{0,f}$ by fixing the $i$th variable
$[x_i,y_i]:=[\alpha,\gamma]$
is reducible.\\
{\rm (c)} $Q_{0,f}$ has an irreducible $(n-2)$-ary retract.
\end{Theorem}
\begin{Proof}
The theorem is a corollary of the properties of the function $f$
discussed in the next section.

(a) By Lemma~\ref{l:0}, the claim follows directly from Theorem~\ref{t:b}(a).

(b) It is straightforward that $Q^i_{[\alpha,\gamma]}=Q_{\alpha,f^{i}_{\alpha}+\gamma}$
where $f^{i}_{\alpha}$ is obtained from $f$
by fixing the $i$th variable $x_i:=\alpha$.
So, by Lemma~\ref{l:0}, the reducibility of $Q^i_{[\alpha,\gamma]}$ is a corollary of
the separability of $f^{i}_{\alpha}$ (Theorem~\ref{t:b}(b)).

Similarly, (c) follows from the fact that fixing two variables we can get
a non-separable subfunction of $f$ (Theorem~\ref{t:b}(c)).
\end{Proof}

\begin{Note}\label{r:red}
An $n$-qua\-si\-group is called \emph{$(i,j)$-re\-du\-ci\-ble} if 
it is $\{i,\ldots,i+j-1\}$-re\-du\-ci\-ble
for some $i\in \{1,\ldots,n\}$ and $j\in \{2,\ldots,n-1\}$ meeting $i+j-1\leq n$.
Clearly, the property of $(i,j)$-re\-du\-ci\-bi\-li\-ty  is stronger
than the permutable reducibility and
is not invariant under changing the argument order;
this property was considered e.\,g. in \cite{Belousov}.
Using an appropriate argument permutation
(more precisely, replacing $f$ by $f'(x_0,x_1,\ldots,x_{n}) \eqdf f(x_0,x_2,\ldots,x_{2n\bmod {(n+1)}})$),
we can strengthen the statement of Theorem~{\rm\ref{t:q}(b)}
getting the $(i,j)$-re\-du\-ci\-ble $(n-1)$-ary retracts.
\end{Note}
\begin{Note}\label{r:order}
Using $Q_{0,f}$ {\rm(}or $Q_{0,f'}$, see Remark~{\rm\ref{r:red}}{\rm)},
it is not difficult to construct an irreducible $n$-qua\-si\-group of order
$4r$ with reducible {\rm(}$(i,j)$-reducible{\rm)} $(n-1)$-ary retracts for any $r>0$: \,
if $(G,*)$ is a commutative group of order $|G|=r\leq\infty$
then the $n$-qua\-si\-group $Q_{f}^{(G,*)}$ {\rm (}and, similarly, its retracts{\rm )}
defined as 
\begin{equation}\label{e:defQG}
Q_{f}^{(G,*)}([w,z]_{1}^{n}) \eqdf [w_1*\ldots*w_n,\,Q_{0,f}(z_{1}^{n})], \qquad w_i \in G,\quad z_i\in Z_2^2
\end{equation}
inherits all the reducibility properties of $Q_{0,f}$ {\rm (}and its retracts{\rm)}.
Indeed, if $Q_{0,f}$ is $A$-re\-du\-ci\-ble then, obviously, $Q_{f}^{(G,*)}$ is $A$-re\-du\-ci\-ble too.
Conversely, let $Q_{f}^{(G,*)}$ be $A$-re\-du\-ci\-ble.
Since the group ${(G,*)}$
is commutative, we can assume without loss of generality that $A=\{\overline{0,k-1}\}$.
Using Lemma~\ref{l:abg}, we can check that
$$
Q_{f}^{(G,*)}([w,z]_{1}^{n}) \equiv
[w_1*\ldots*w_n,\,
Q_{0,f}(z_{1}^{k-1},q^{-1}(Q_{0,f}(0^{k-1},z_{k}^{n})),0^{n-k})]
$$
with $q(z) \eqdf Q_{0,f}(0^{k-1},z,0^{n-k})$.
Comparing with (\ref{e:defQG}) gives a reduction of $Q_{0,f}$.
%
\end{Note}

\section{Properties of the partial Boolean function $\lowercase{f}$} \label{s:b}

In this section we prove the key theorem of the paper:
\begin{Theorem}\label{t:b}
Let $n\geq 4$ be even and the partial Boolean function $f:E^{n+1}_{0}\to E$ be represented
by the following polynomial$:$
\begin{equation}\label{e:f}
f(x_{0}^{n}) \eqdf \sum_{i=0}^{n}\sum_{j=1}^{\lfloor n/4 \rfloor} x_i x_{i+j}
\end{equation}
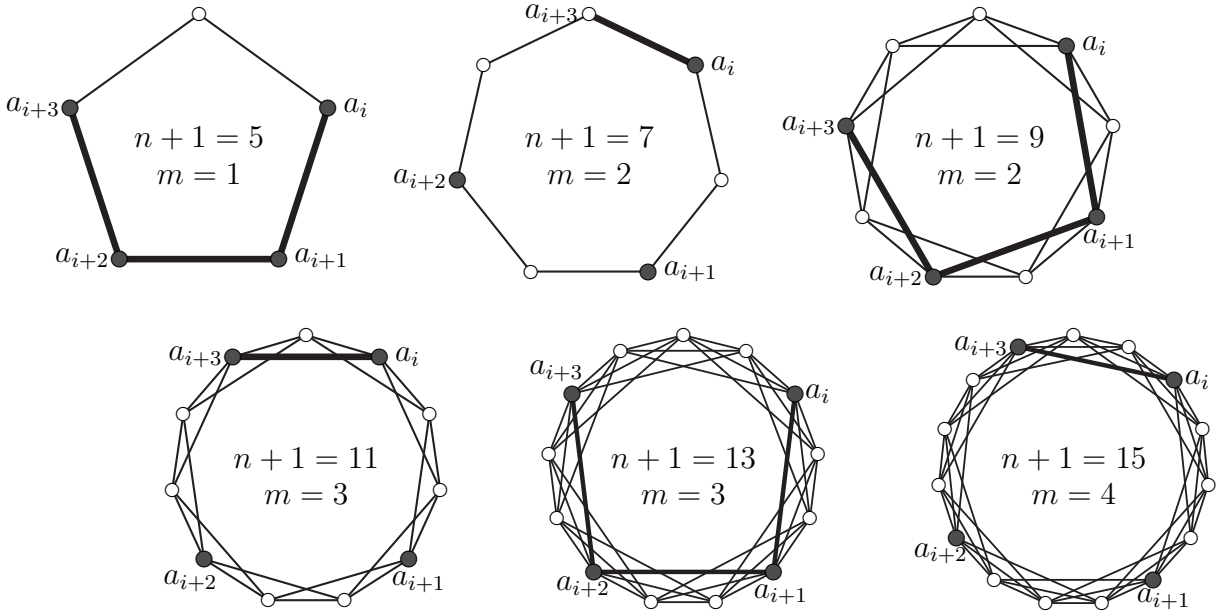
\begin{figure}
\begin{center}
 \unitlength 1.8mm 
\def\EMmv#1#2{\put(#1,#2){\special{em:moveto}}}
\def\EMln#1#2{\put(#1,#2){\special{em:lineto}}}
\def\EMnorm{\special{em:linewidth 0.8pt}} 
\def\EMdick{\special{em:linewidth 2.5pt}} 
\def\PUTcr#1#2#3{\put(#1,#2){\circle{#3}}}
\def\PUTdt#1#2#3{\put(#1,#2){\circle*{#3}}}
\def\PUTvx#1#2#3{\color{white}\put(#1,#2){\circle*{1}}\color{black}\put(#1,#2){\circle{1}}}
\def\PUTvv#1#2#3{\color[rgb]{0.3,0.3,0.3}\put(#1,#2){\circle*{1.2}}\color{black}\put(#1,#2){\circle{1.2}}}
\begin{picture}(25,21)(12,0)
\def\xxa{10,000} \def\yya{20,000}
\def\xxb{19,511} \def\yyb{13,090}
\def\xxc{15,878} \def\yyc{01,910}
\def\xxd{04,122} \def\yyd{01,910}
\def\xxe{00,489} \def\yye{13,090}
\EMnorm
\EMmv{\xxa}{\yya}
\EMln{\xxb}{\yyb}
\EMdick
\EMln{\xxc}{\yyc}
\EMln{\xxd}{\yyd}
\EMln{\xxe}{\yye}
\EMnorm
\EMln{\xxa}{\yya}
\EMdick
\PUTvx{\xxa}{\yya}{}
\PUTvv{\xxb}{\yyb}{}
\PUTvv{\xxc}{\yyc}{}
\PUTvv{\xxd}{\yyd}{}
\PUTvv{\xxe}{\yye}{}
\put(10,10){\makebox(0,0)[bc]{$n+1=5$}}
\put(10, 9){\makebox(0,0)[tc]{$m=1$}}
\put(\xxb,\yyb){\makebox(0,0)[cl]{$\ \,a_i$}}
\put(\xxc,\yyc){\raisebox{0mm}{\makebox(0,0)[cl]{$\ \,a_{i+1}$}}}
\put(\xxd,\yyd){\makebox(0,0)[cr]{$a_{i+2}\ $}}
\put(\xxe,\yye){\makebox(0,0)[cr]{$a_{i+3}\ $}}
\end{picture}
\begin{picture}(25,21)(12,0)
\def\xxa{10,000} \def\yya{20,000}
\def\xxb{17,818} \def\yyb{16,235}
\def\xxc{19,749} \def\yyc{07,775}
\def\xxd{14,339} \def\yyd{00,990}
\def\xxe{05,661} \def\yye{00,990}
\def\xxf{00,251} \def\yyf{07,775}
\def\xxg{02,182} \def\yyg{16,235}
\EMnorm
\EMdick
\EMmv{\xxa}{\yya}
\EMln{\xxb}{\yyb}
\EMnorm
\EMln{\xxc}{\yyc}
\EMln{\xxd}{\yyd}
\EMln{\xxe}{\yye}
\EMln{\xxf}{\yyf}
\EMln{\xxg}{\yyg}
\EMln{\xxa}{\yya}
\EMdick
\PUTvx{\xxa}{\yya}{}
\PUTvv{\xxb}{\yyb}{}
\PUTvx{\xxc}{\yyc}{}
\PUTvv{\xxd}{\yyd}{}
\PUTvx{\xxe}{\yye}{}
\PUTvv{\xxf}{\yyf}{}
\PUTvx{\xxg}{\yyg}{}
\put(10,10){\makebox(0,0)[bc]{$n+1=7$}}
\put(10, 9){\makebox(0,0)[tc]{$m=2$}}
\put(\xxb,\yyb){\makebox(0,0)[cl]{$\ \,a_i$}}
\put(\xxd,\yyd){\raisebox{0mm}{\makebox(0,0)[cl]{$\ \,a_{i+1}$}}}
\put(\xxf,\yyf){\makebox(0,0)[cr]{$a_{i+2}\ $}}
\put(\xxa,\yya){\makebox(0,0)[cr]{$a_{i+3}\ $}}
\end{picture}
\begin{picture}(25,21)(12,0)
\def\xxa{10,000} \def\yya{20,000}
\def\xxb{16,428} \def\yyb{17,660}
\def\xxc{19,848} \def\yyc{11,736}
\def\xxd{18,660} \def\yyd{05,000}
\def\xxe{13,420} \def\yye{00,603}
\def\xxf{06,580} \def\yyf{00,603}
\def\xxg{01,340} \def\yyg{05,000}
\def\xxh{00,152} \def\yyh{11,736}
\def\xxi{03,572} \def\yyi{17,660}
\EMnorm
\EMmv{\xxa}{\yya}
\EMln{\xxb}{\yyb}
\EMln{\xxc}{\yyc}
\EMln{\xxd}{\yyd}
\EMln{\xxe}{\yye}
\EMln{\xxf}{\yyf}
\EMln{\xxg}{\yyg}
\EMln{\xxh}{\yyh}
\EMln{\xxi}{\yyi}
\EMln{\xxa}{\yya}
\EMln{\xxc}{\yyc}
\EMln{\xxe}{\yye}
\EMln{\xxg}{\yyg}
\EMln{\xxi}{\yyi}
\EMln{\xxb}{\yyb}
\EMdick
\EMln{\xxd}{\yyd}
\EMln{\xxf}{\yyf}
\EMln{\xxh}{\yyh}
\EMnorm
\EMln{\xxa}{\yya}
\EMdick
\PUTvx{\xxa}{\yya}{}
\PUTvv{\xxb}{\yyb}{}
\PUTvx{\xxc}{\yyc}{}
\PUTvv{\xxd}{\yyd}{}
\PUTvx{\xxe}{\yye}{}
\PUTvv{\xxf}{\yyf}{}
\PUTvx{\xxg}{\yyg}{}
\PUTvv{\xxh}{\yyh}{}
\PUTvx{\xxi}{\yyi}{}
\put(10,10){\makebox(0,0)[bc]{$n+1=9$}}
\put(10, 9){\makebox(0,0)[tc]{$m=2$}}
\put(\xxb,\yyb){\makebox(0,0)[cl]{$\ \,a_i$}}
\put(\xxd,\yyd){\raisebox{-3mm}{\makebox(0,0)[cc]{$\ \ \,a_{i+1}$}}}
\put(\xxf,\yyf){\makebox(0,0)[cr]{$a_{i+2}\ $}}
\put(\xxh,\yyh){\makebox(0,0)[cr]{$a_{i+3}\ $}}
\end{picture}
\\[2.5ex]
\begin{picture}(25,21)(12.5,0)
\def\xxa{10,000} \def\yya{20,000}
\def\xxb{15,406} \def\yyb{18,413}
\def\xxc{19,096} \def\yyc{14,154}
\def\xxd{19,898} \def\yyd{08,577}
\def\xxe{17,557} \def\yye{03,451}
\def\xxf{12,817} \def\yyf{00,405}
\def\xxg{07,183} \def\yyg{00,405}
\def\xxh{02,443} \def\yyh{03,451}
\def\xxi{00,102} \def\yyi{08,577}
\def\xxj{00,904} \def\yyj{14,154}
\def\xxk{04,594} \def\yyk{18,413}
\EMnorm
\EMmv{\xxa}{\yya}
\EMln{\xxb}{\yyb}
\EMln{\xxc}{\yyc}
\EMln{\xxd}{\yyd}
\EMln{\xxe}{\yye}
\EMln{\xxf}{\yyf}
\EMln{\xxg}{\yyg}
\EMln{\xxh}{\yyh}
\EMln{\xxi}{\yyi}
\EMln{\xxj}{\yyj}
\EMln{\xxk}{\yyk}
\EMln{\xxa}{\yya}
\EMln{\xxc}{\yyc}
\EMln{\xxe}{\yye}
\EMln{\xxg}{\yyg}
\EMln{\xxi}{\yyi}
\EMln{\xxk}{\yyk}
\EMdick
\EMln{\xxb}{\yyb}
\EMnorm
\EMln{\xxd}{\yyd}
\EMln{\xxf}{\yyf}
\EMln{\xxh}{\yyh}
\EMln{\xxj}{\yyj}
\EMln{\xxa}{\yya}
\EMdick
\PUTvx{\xxa}{\yya}{}
\PUTvv{\xxb}{\yyb}{}
\PUTvx{\xxc}{\yyc}{}
\PUTvx{\xxd}{\yyd}{}
\PUTvv{\xxe}{\yye}{}
\PUTvx{\xxf}{\yyf}{}
\PUTvx{\xxg}{\yyg}{}
\PUTvv{\xxh}{\yyh}{}
\PUTvx{\xxi}{\yyi}{}
\PUTvx{\xxj}{\yyj}{}
\PUTvv{\xxk}{\yyk}{}
\put(10,10){\makebox(0,0)[bc]{$n+1=11$}}
\put(10, 9){\makebox(0,0)[tc]{$m=3$}}
\put(\xxb,\yyb){\makebox(0,0)[cl]{$\ \,a_i$}}
\put(\xxe,\yye){\raisebox{-3mm}{\makebox(0,0)[cc]{$\ \ a_{i+1}$}}}
\put(\xxh,\yyh){\raisebox{-1.5mm}{\makebox(0,0)[tr]{$a_{i+2}$\hspace{-1ex}}}}
\put(\xxk,\yyk){\makebox(0,0)[cr]{$a_{i+3}\ $}}
\end{picture}
\def\EMnorm{\special{em:linewidth 0.7pt}} 
\def\EMdick{\special{em:linewidth 1.7pt}} 
\begin{picture}(25,21)(15,0)
\def\xxa{10,000} \def\yya{20,000}
\def\xxb{14,647} \def\yyb{18,855}
\def\xxc{18,230} \def\yyc{15,681}
\def\xxd{19,927} \def\yyd{11,205}
\def\xxe{19,350} \def\yye{06,454}
\def\xxf{16,631} \def\yyf{02,515}
\def\xxg{12,393} \def\yyg{00,291}
\def\xxh{07,607} \def\yyh{00,291}
\def\xxi{03,369} \def\yyi{02,515}
\def\xxj{00,650} \def\yyj{06,454}
\def\xxk{00,073} \def\yyk{11,205}
\def\xxl{01,770} \def\yyl{15,681}
\def\xxm{05,353} \def\yym{18,855}
\EMnorm
\EMmv{\xxa}{\yya}
\EMln{\xxb}{\yyb}
\EMln{\xxc}{\yyc}
\EMln{\xxd}{\yyd}
\EMln{\xxe}{\yye}
\EMln{\xxf}{\yyf}
\EMln{\xxg}{\yyg}
\EMln{\xxh}{\yyh}
\EMln{\xxi}{\yyi}
\EMln{\xxj}{\yyj}
\EMln{\xxk}{\yyk}
\EMln{\xxl}{\yyl}
\EMln{\xxm}{\yym}
\EMln{\xxa}{\yya}
\EMln{\xxc}{\yyc}
\EMln{\xxe}{\yye}
\EMln{\xxg}{\yyg}
\EMln{\xxi}{\yyi}
\EMln{\xxk}{\yyk}
\EMln{\xxm}{\yym}
\EMln{\xxb}{\yyb}
\EMln{\xxd}{\yyd}
\EMln{\xxf}{\yyf}
\EMln{\xxh}{\yyh}
\EMln{\xxj}{\yyj}
\EMln{\xxl}{\yyl}
\EMln{\xxa}{\yya}
\EMln{\xxd}{\yyd}
\EMln{\xxg}{\yyg}
\EMln{\xxj}{\yyj}
\EMln{\xxm}{\yym}
\EMln{\xxc}{\yyc}
\EMdick
\EMln{\xxf}{\yyf}
\EMln{\xxi}{\yyi}
\EMln{\xxl}{\yyl}
\EMnorm
\EMln{\xxb}{\yyb}
\EMln{\xxe}{\yye}
\EMln{\xxh}{\yyh}
\EMln{\xxk}{\yyk}
\EMln{\xxa}{\yya}
\EMdick
\PUTvx{\xxa}{\yya}{}
\PUTvx{\xxb}{\yyb}{}
\PUTvv{\xxc}{\yyc}{}
\PUTvx{\xxd}{\yyd}{}
\PUTvx{\xxe}{\yye}{}
\PUTvv{\xxf}{\yyf}{}
\PUTvx{\xxg}{\yyg}{}
\PUTvx{\xxh}{\yyh}{}
\PUTvv{\xxi}{\yyi}{}
\PUTvx{\xxj}{\yyj}{}
\PUTvx{\xxk}{\yyk}{}
\PUTvv{\xxl}{\yyl}{}
\PUTvx{\xxm}{\yym}{}
\put(10,10){\makebox(0,0)[bc]{$n+1=13$}}
\put(10, 9){\makebox(0,0)[tc]{$m=3$}}
\put(\xxc,\yyc){\makebox(0,0)[cl]{$\ a_i$}}
\put(\xxf,\yyf){\raisebox{-3mm}{\makebox(0,0)[cc]{$\ \,a_{i+1}$}}}
\put(\xxi,\yyi){\raisebox{-2mm}{\makebox(0,0)[cc]{$a_{i+2}\ \ $}}}
\put(\xxl,\yyl){\raisebox{3mm}{\makebox(0,0)[cc]{$a_{i+3}\ \ \ $}}}
\end{picture}
\begin{picture}(25,21)(15,0)
\def\xxa{10,000} \def\yya{20,000}
\def\xxb{14,067} \def\yyb{19,135}
\def\xxc{17,431} \def\yyc{16,691}
\def\xxd{19,511} \def\yyd{13,090}
\def\xxe{19,945} \def\yye{08,955}
\def\xxf{18,660} \def\yyf{05,000}
\def\xxg{15,878} \def\yyg{01,910}
\def\xxh{12,079} \def\yyh{00,219}
\def\xxi{07,921} \def\yyi{00,219}
\def\xxj{04,122} \def\yyj{01,910}
\def\xxk{01,340} \def\yyk{05,000}
\def\xxl{00,055} \def\yyl{08,955}
\def\xxm{00,489} \def\yym{13,090}
\def\xxn{02,569} \def\yyn{16,691}
\def\xxo{05,933} \def\yyo{19,135}
\EMnorm
\EMmv{\xxa}{\yya}
\EMln{\xxb}{\yyb}
\EMln{\xxc}{\yyc}
\EMln{\xxd}{\yyd}
\EMln{\xxe}{\yye}
\EMln{\xxf}{\yyf}
\EMln{\xxg}{\yyg}
\EMln{\xxh}{\yyh}
\EMln{\xxi}{\yyi}
\EMln{\xxj}{\yyj}
\EMln{\xxk}{\yyk}
\EMln{\xxl}{\yyl}
\EMln{\xxm}{\yym}
\EMln{\xxn}{\yyn}
\EMln{\xxo}{\yyo}
\EMln{\xxa}{\yya}
\EMln{\xxc}{\yyc}
\EMln{\xxe}{\yye}
\EMln{\xxg}{\yyg}
\EMln{\xxi}{\yyi}
\EMln{\xxk}{\yyk}
\EMln{\xxm}{\yym}
\EMln{\xxo}{\yyo}
\EMln{\xxb}{\yyb}
\EMln{\xxd}{\yyd}
\EMln{\xxf}{\yyf}
\EMln{\xxh}{\yyh}
\EMln{\xxj}{\yyj}
\EMln{\xxl}{\yyl}
\EMln{\xxn}{\yyn}
\EMln{\xxa}{\yya}
\EMln{\xxd}{\yyd}
\EMln{\xxg}{\yyg}
\EMln{\xxj}{\yyj}
\EMln{\xxm}{\yym}
\EMln{\xxa}{\yya}
\EMmv{\xxb}{\yyb}
\EMln{\xxe}{\yye}
\EMln{\xxh}{\yyh}
\EMln{\xxk}{\yyk}
\EMln{\xxn}{\yyn}
\EMln{\xxb}{\yyb}
\EMmv{\xxc}{\yyc}
\EMln{\xxf}{\yyf}
\EMln{\xxi}{\yyi}
\EMln{\xxl}{\yyl}
\EMln{\xxo}{\yyo}
\EMdick
\EMln{\xxc}{\yyc}
\EMnorm
\EMdick
\PUTvx{\xxa}{\yya}{}
\PUTvx{\xxb}{\yyb}{}
\PUTvv{\xxc}{\yyc}{}
\PUTvx{\xxd}{\yyd}{}
\PUTvx{\xxe}{\yye}{}
\PUTvx{\xxf}{\yyf}{}
\PUTvv{\xxg}{\yyg}{}
\PUTvx{\xxh}{\yyh}{}
\PUTvx{\xxi}{\yyi}{}
\PUTvx{\xxj}{\yyj}{}
\PUTvv{\xxk}{\yyk}{}
\PUTvx{\xxl}{\yyl}{}
\PUTvx{\xxm}{\yym}{}
\PUTvx{\xxn}{\yyn}{}
\PUTvv{\xxo}{\yyo}{}
\put(10,10){\makebox(0,0)[bc]{$n+1=15$}}
\put(10, 9){\makebox(0,0)[tc]{$m=4$}}
\put(\xxc,\yyc){\makebox(0,0)[cl]{$\ a_i$}}
\put(\xxg,\yyg){\raisebox{-2.5mm}{\makebox(0,0)[cc]{$\ \ a_{i+1}$}}}
\put(\xxk,\yyk){\raisebox{-2mm}{\makebox(0,0)[cc]{$a_{i+2}\ \ \ $}}}
\put(\xxo,\yyo){\raisebox{0mm}{\makebox(0,0)[cr]{$a_{i+3}\ $}}}
\end{picture}
 \caption{\label{f:} It is natural to represent a square-free
 (i.\,e., without monomials of type $x_i^2$) quadratic form over $Z_2$
 by the graph
 whose $i$th and $j$th vertices are connected
 if and only if the form contains the monomial $x_i x_j$.
 The figure presents the graph corresponding to the form (\ref{e:f}) with $n=4$, $6$, $8$, $10$, $12$, and $14$.
  }
\end{center}
\end{figure}
$($see Fig.~\ref{f:}$)$. Put $m \eqdf \lfloor (n+2)/4 \rfloor$. Then\\
{\rm (a)} The partial Boolean function $f$ is not separable.\\
{\rm (b)} For all $i\in \{\overline{0,n}\}$ and $\alpha\in \{0,1\}$ the subfunction
$f^i_{\alpha}:E^{n}_{\alpha}\to E$
obtained from $f(x_{0}^{n})$ by fixing 
$x_i:=\alpha$
is $\{{i+m},{i-m}\}$-se\-pa\-rable
$($here and in what follows for subfunctions we leave the same numeration of variables
as for the original function$)$.\\
{\rm (c)} For all $i\in \{\overline{0,n}\}$ and $\alpha,\beta\in \{0,1\}$
the subfunction $g^i_{\alpha,\beta}:E^{n-1}_{\alpha+\beta}\to E$
obtained from $f(x_{0}^{n})$ by fixing 
$x_i:=\alpha$, $x_{i+m}:=\beta$ is not separable.
\end{Theorem}
\begin{Proof}
(a) Let $A$ be an arbitrary subset of $\{\overline{0,n}\}$ such that $2\leq|A|\leq n-1$,
and let $B\eqdf \{\overline{0,n}\}\backslash A$.
We will show that $f$ is not $A$-se\-pa\-rable, using the two following simple facts:

\begin{Lemma}\label{p:sub}
Assume a partial Boolean function $f:E^{n+1}_0\to \{0,1\}$ is $A$-se\-pa\-rable.
Then each $($partial$)$ subfunction $f'$ obtained from $f(x_{0}^{n})$ by fixing some variables
$x_{v_1},\ldots,x_{v_k}$ is
$A'$-se\-pa\-rable with $A'\eqdf A\backslash \{v_{1}^{k}\}$.
\end{Lemma}
\begin{Lemma}\label{p:0123}
Let $\gamma_{01},\gamma_{02},\gamma_{03},\gamma_{12},\gamma_{13},\gamma_{23}\in \{0,1\}$.
A partial Boolean function 
$$
h(x_0,x_1,x_2,x_3)\eqdf\gamma_{01}x_0 x_1+\gamma_{02}x_0 x_2+\gamma_{03}x_0 x_3
+\gamma_{12}x_1 x_2+\gamma_{13}x_1 x_3+\gamma_{23}x_2 x_3:
\qquad
$$
$E^4_0\to \{0,1\}$ is $\{0,1\}$-se\-pa\-rable
only if $\gamma_{02}+\gamma_{03}+\gamma_{12}+\gamma_{13}=0$.
\end{Lemma}
\noindent (Lemma~\ref{p:sub} is straightforward from the definition.
Proof of Lemma~\ref{p:0123}:
From the $\{0,1\}$-se\-pa\-ra\-bi\-li\-ty of $h$ we derive
$h(0,0,\,0,0)+h(1,1,\,1,1)=h(1,1,\,0,0)+h(0,0,\,1,1)$.
Substituting the definition of $h$, we get
$\gamma_{02}+\gamma_{03}+\gamma_{12}+\gamma_{13}=0$.)


Consider the cyclic sequence $a_i=i\cdot m \bmod (n+1)$, $i=0,\ldots, n$.
Since $n+1=4m\pm 1$, we see that $m$ and $n+1$ are relatively prime, and $\{a_{0}^{n}\}=\{\overline{0,n}\}$.
At least one of the following holds (recall that indices are calculated modulo~$n+1$):

1) $a_i,a_{i+1}\in A$, $a_{i+2},a_{i+3}\in B$
or $a_i,a_{i+1}\in B$, $a_{i+2},a_{i+3}\in A$
for some $i$.
Assigning zeroes to all variables of $f(x_{0}^{n})$ except
$x_{a_{i}},x_{a_{i+1}},x_{a_{i+2}},x_{a_{i+3}}$ we get the partial Boolean function
$$
f'(x_{a_{i}},x_{a_{i+1}},x_{a_{i+2}},x_{a_{i+3}}) \equiv
\cases{
x_{a_{i}}x_{a_{i+1}}+x_{a_{i+1}}x_{a_{i+2}}+x_{a_{i+2}}x_{a_{i+3}},
                      & \hspace{-0em} if $n\equiv 0 \bmod 4$,\cr
x_{a_{i}}x_{a_{i+3}}, & \hspace{-0em} if $n\equiv 2 \bmod 4$
}
$$
(see Fig.~\ref{f:}, the dark nodes),
which is not $\{a_{i},a_{i+1}\}$-se\-pa\-rable, by Lemma~\ref{p:0123}.
Therefore $f$ is not $A$-se\-pa\-rable, by Lemma~\ref{p:sub}.

2) $a_i,a_{i+2}\in A$, $a_{i+1}\in B$
or $a_i,a_{i+2}\in B$, $a_{i+1}\in A$
for some $i$.
Without loss of generality assume $0\in A$, $m\in B$, $2m\in A$.
Note that the polynomial (\ref{e:f}) contains exactly one of monomials
$x_{0}x_{b}$, $x_{2m}x_{b}$ for each $b\neq 0,m,2m$.
Take $b\in B\backslash\{ m \}$. 
Assigning zeroes to all variables of $f(x_{0}^{n})$ except
$x_{0},x_{m},x_{2m},x_{b}$ we get the partial Boolean function
$$
f''(x_{0},x_{2m},x_{m},x_{b})
 \equiv  \cases{
x_{0}x_{m}+x_{m}x_{2m}+\alpha x_{0}x_{b}+\beta x_{m}x_{b}+{\bar \alpha} x_{2m}x_{b},
                                                             & \hspace{-0em} if $n\equiv 0 \bmod 4$,\cr
\alpha x_{0}x_{b}+\beta x_{m}x_{b}+{\bar \alpha} x_{2m}x_{b},& \hspace{-0em} if $n\equiv 2 \bmod 4$
}
$$
with $\alpha, \beta\in \{0,1\}$, ${\bar \alpha}\eqdf 1-\alpha$.
In any case, $f''(x_{0},x_{m},x_{2m},x_{b})$ is not
$\{0,2m\}$-se\-pa\-rable, by Lemma~\ref{p:0123}.
It follows that $f$ is not $A$-se\-pa\-rable, by Lemma~\ref{p:sub}.

(b) Without loss of generality we assume $i=0$.
Put
$$\tilde x_k\eqdf |x_{k-\lfloor n/4 \rfloor}^{k-1}| + |x_{k+1}^{k+\lfloor n/4 \rfloor}|
=|x_{k-\lfloor n/4 \rfloor}^{k+\lfloor n/4 \rfloor}|+x_k.$$
Note that $m+\lfloor n/4 \rfloor=n/2$, and $m-\lfloor n/4 \rfloor$ is $0$ or $1$;
in both cases,
$$
    |x_{0}^{n}| \equiv \big(\tilde x_{m}+x_{m}+\tilde x_{-m}+x_{-m}+x_{0}\big).
$$
Since $|x_{0}^{n}|$ equals zero everywhere on $E^{n+1}_0$,
we can represent $f$ as follows:
\begin{eqnarray*}
f(x_{0}^{n})
&\equiv & \sum_{i=0}^{n}\sum_{j=1}^{\lfloor n/4 \rfloor} x_i x_{i+j}
+ \big(\tilde x_{m}+x_{m}+\tilde x_{-m}+x_{-m}+x_{0}\big)
\big(\tilde x_{m} + x_{-m}\big)\\
&\equiv & \sum_{i=0}^{n}\sum_{j=1}^{\lfloor n/4 \rfloor} x_i x_{i+j}
+x_{m}\tilde x_{m}
+x_{-m}\tilde x_{-m}
+(x_{m}+x_{-m}+x_{0})x_{-m}+{S}
\end{eqnarray*}
where $S$ does not depend on $x_{m}$ and $x_{-m}$.
It is easy to see that this representation does not contain any monomial
$x_k x_{k'}$ with $k\in \{-m,m\}$, $k'\not\in \{0,-m,m\}$.
This means that after fixing $x_0$ we obtain a $\{-m,m\}$-se\-pa\-rable partial Boolean function.

(c) Without loss of generality assume $i=0$.
Let $A$ be an arbitrary subset of $\{\overline{1,m-1},\overline{m+1,n}\}$ such that $2\leq|A|\leq n-2$;
let $B\eqdf \{\overline{1,m-1},\overline{m+1,n}\}\backslash A$.
If the sequence $a_i$, $i=\overline{0,n}$ is defined as in (a) then
either 1) or 2) holds or

3) $A=\{a_2, a_{n}\} = \{2m, -m\}$ or $B= \{2m, -m\}$
(recall that the numbers $a_0=0$ and $a_1=m$ correspond to the fixed variables).
As in the cases 1) and 2),
assigning zeroes to all variables of
$g^0_{\alpha,\beta}(x_{1}^{m-1},x_{m+1}^{n})=f(\alpha,x_{1}^{m-1},\beta,x_{m+1}^{n})$ except
$x_{2m},x_{-m},x_{1},x_{n}$, we find that
$g^0_{\alpha,\beta}$ is not $A$-se\-pa\-rable by Lemmas~\ref{p:sub} and~\ref{p:0123}.
\end{Proof}
In the proof of the part (b) we exploit the fact that after removing a vertex, say $0$,
in the corresponding graph (see Fig.~\ref{f:}) the remaining vertex set will be the
disjoint union of the two vertices $m$ and $-m$ and their neighborhoods.
This partly explains why our construction does not work in the case of even $n+1$.
In the following remark we compare our results with the situation with (total) Boolean functions.
\begin{Note}
Say that a Boolean function $\mu(x_1,\ldots,x_n):E^n\to\{0,1\}$ is \emph{separable} if 
it is $A$-se\-pa\-rable for some
$A\subset \{\overline{1,n}\}$ where $1\leq|A|\leq n-1$ and $A$-se\-pa\-rabi\-li\-ty means
the same as for partial Boolean functions.
Then (*) every non-se\-pa\-rable $n$-ary Boolean function $\mu$ has
a non-se\-pa\-rable $(n-1)$-ary subfunction obtained from $\mu$ by fixing some variable.
(Assume the contrary;
consider a maximal non-se\-pa\-rable $k$-ary subfunction $\mu'$;
and prove that $\mu=\mu'+\mu''$ for some $(n-k)$-ary $\mu''$
where the free variables in $\mu'$ and $\mu''$ do not intersect).
Our investigation shows that the situation with the partial Boolean functions on $E^{n+1}_0$
is 
more complex; the statement like (*) fails for even $n$ and holds for $n=5$ and $n=7$.
Question: does it hold for every odd $n$?
\end{Note}

\section{Remark. Switching separability of graphs}\label{s:g}

As noted in the comments on Fig.~\ref{f:},
each square-free quadratic form $p(x_{0}^{n})$ over $Z_2$
can be represented by the graph with
$n+1$ vertices $\{0,\ldots,n\}$ such that vertices $i$ and $j$ are adjacent
if and only if $p(x_{0}^{n})$ contains the monomial $x_i x_j$.
In this section we define the concept of graph switching separability that corresponds
to the separability of the corresponding quadratic polynomial considered as a partial
Boolean function $E^{n+1}_0 \to \{0,1\}$.

 We first define a graph transformation,
which is known as a \emph{graph switching} or \emph{Seidel switching}.
The result of \emph{switching} a set $U\subseteq V$ in a graph
$G=(V,E)$ is defined as the graph with the same vertex set $V$ and the edge set
$E\vartriangle E_{U,V\setminus U}$ where
$E_{U,V\setminus U}\eqdf \{\{u,v\}\,|\,u\in U, v\in V\setminus U\}$.
We say that the graph $G=(V,E)$ is \emph{switching-separable} if 
$V=V_1\cup V_2$ where $|V_1|\geq 2$, $|V_2|\geq 2$, $V_1\cap V_2=\emptyset$, and
for some $U\subseteq V$ switching $U$ in $G$ gives a graph with no edges between
$V_1$ and $V_2$. Clearly, if a graph is switching-separable then
all its switchings are switching-separable. The class of all switchings
of a graph is known as a \emph{switchings class} and is equivalent to a \emph{two-graph},
see e.\,g. \cite{Spence:2graphs}.
From Theorem~\ref{t:b} and the computer search observed in the Introduction,
we can derive the following:
\begin{Corollary}
\label{c:g}
For every odd $|V|\geq 5$ there exists a non switching-separable graph
$G=(V,E)$ such that every subgraph generated by $|V|-1$ vertices is switching-separable. If
$|V|=6$ or $|V|=8$ then such graphs do not exist.
\end{Corollary}



\providecommand\href[2]{#2}\providecommand\url[1]{\href{#1}{#1}}

\end{document}